# Optimization as a design strategy.
# Considerations based on building simulation-assisted experiments about problem decomposition

Gian Luca Brunetti[1]

**Abstract.** In this article the most fundamental decomposition-based optimization method - block coordinate search, based on the sequential decomposition of problems in subproblems - and building performance simulation programs are used to reason about a building design process at micro-urban scale and strategies are defined to make the search more efficient. Cyclic overlapping block coordinate search is here considered in its double nature of optimization method and surrogate model (and metaphore) of a sequential design process. Heuristic indicators apt to support the design of search structures suited to that method are developed from building-simulation-assisted computational experiments, aimed to choose the form and position of a small building in a plot. Those indicators link the sharing of structure between subspaces ("commonality") to recursive recombination, measured as freshness of the search wake and novelty of the search moves. The aim of these indicators is to measure the relative effectiveness of decomposition-based design moves and create efficient block searches. Implications of a possible use of these indicators in genetic algorithms are also highlighted.

## 1. Introduction

Building and urban design are today increasingly supported by optimization methods. In this paper criteria will be defined for adding efficiency to sequential optimization searches consisting in the decomposition of problems in (possibly overlapping) subproblems; and elements of reasoning will be sought to ascertain if those criteria may be used in population-based optimization methods, which are stemmed from a historical split from enumerative techniques.

---

[1] Department of Architecture and Urban Studies (DAStU), Polytechnic of Milan, gianluca.brunetti@polimi.it.

Historical splits also concern building design. In building design, an awareness about the design process commensurate to the available contemporary knowledge ceased to be attainable by individuals when the scientific progress made during the Age of Enlightenment prompted the separation of competences between building as a technique, which may be considered more typical of engineers, and building as a craft, which may be considered more typical of architects (Kostof, 1977; Cowan, 1978, 1985). Before that separation was accomplished, the capability of a person sufficed to master the qualitative and quantitative knowledge at one's reach about building design. But after that, the quantity of available knowledge soon soared far from the reach of any single person. A consequence of that separation was that engineers and architects began to design different types of buildings. A more profound one was that it prompted the consolidation of two distinct paths to design. One, mainly adopted by engineers, was principally based on calculations and had the advantage of being more exact, making possible to foresee the results of new, non-experimented solutions. The other, mainly adopted by architects, was based on the patchwork re-use of tested examples, somewhat resembling the strategies of savage thinking (Lévi-Strauss, 1962), and had the advantage of being fast, and therefore suitable to lateral thinking, but the disadvantage of discouraging the experimentation of new solutions. The two types of approach coexist to this day in real-world engineering and architectural practice and a synthesis between them has until now been left to the sensibility of professionals and is mostly acquired through example and trial, very much like a medieval, pre-scientific craft.

A flexibility of assessment which may bring to a new unity the "design by calculation" and the "design by example" approaches may be attained in the future by the convergence of BIM programs and simulation programs. That would be an important milestone, since it would bring anew at disposal of individuals a design vision adequate to the knowledge of their times, allowing them to profit from the available scientific knowledge even without being in full command of it. The achievement of that flexibility would be of the highest importance for the quality of design, because the tools (technologies, conventions and languages included) we rely on shape substantially our way of thinking (Leroi-Gourhan, 1964-65; Ong, 1982), and building design is nowadays increasingly relying on simulation tools.

For the attainment of that objective, the definition of efficient strategies for searching into design spaces is important, but probably yet to complete. Reviews and assessments (Evins, 2013; Attia et al, 2013) show that simulation-assisted optimization strategies in building design are today mostly used as a means of fine-tuning, characterized by clear goals and few variables. A reason for this is that managing parametric explorations flexibly is currently difficult at implementation level. Another one is that optimization tools are today prevalently used as black boxes for approaching whole design spaces without preliminarily decomposing them. But this avoidance of decomposition may not be beneficial: decomposition can constitute a valuable contribution of experts to computer-driven optimization enquiries. Furthermore, it may be useful as an algorithm-driven search strategy. In many research areas, like genetic programming (Koza, 1992), convolutional neural networks (Fukushima, 1980; LeCun, 1995) and cognitive science (Anderson et al, 2011), the hierarchical organization of overlapping blocks of information is today emerging as a key strategy towards efficient intelligence.

**2. Background on optimization**
**2.1. Coexistence of optimization and decomposition**

Optimization is so omnipresent in every aspect of existence that its importance is difficult to overstate. Life itself after Darwin can be seen as a recursive process aimed to maximize the likelihood of itself. And optimization plays a crucial role in building and urban design as well. Optimization processes indeed in essence search processes, like design.

The importance of decomposition is also difficult to overstate. Decomposition is a means of processing information with less processing power; a means of reducing search size without reducing search space. And it is at the basis of both low-level neural processes (Hubel et al, 1959) and high-level human thinking, being at the core of problem reduction; which is a condition for the effectiveness of low-level, short-time-span human reasoning (Miller, 1956) and a key component of concept abstraction (Anderson, 2002).

Problem decomposition in last decades has even made its appearance in information-based theories in life sciences. After Dawkins (1976), the possibility has been introduced that the subjects of life may not be whole organisms, but blocks of genetic code struggling for permanence. In such a view, optimization and decomposition are one.

## 2.2. Objectives and globality

Optimisation is relative to objectives. It can be single-objective or multi-objective. Real-life design problems are usually multi-objective. The prevalent strategy to manage multi-objective problems is to transform them in single-objective ones by compounding all the objective functions into one in which they are weighted for importance or relevance for search performance; which can be problematic when objectives are conflicting (Diakaki et al, 2008). Approaches like *physical programming* have gone beyond linear weighting of objective functions, through the use of conditions and non-linear mathematical relations and constraints (Messac, 1996, 2000, Messac et al, 2000, Marler et al, 2004). Population-based search methods, based on populations of coexisting solutions, can instead pursue different objectives in parallel, so to define Pareto non-dominated solutions without requiring an explicit weighting of objective functions (Goldberg, 1989).

Optimization strategies can be classified as local and global ones. Local ones search for local optima and are apt to efficiently converge towards them, at the risk of getting trapped into them. Global search strategies are less efficient in convergence toward local optima, but less subject to get trapped into local optima when dealing with multiconvex or discontinuous problems.

## 2.3. Derivative-free methods

The increasing reliance on simulation programs has determined in the last decades a rise of interest in direct search methods (Wright, 1995), due to the fact that the results produced through them are mainly discontinuous and non-derivable (Kolda et al, 2003, Lewis et al, 2000, Wetter, 2005). Direct search contains the type of pattern search (Hooke et al, 1961, Torczon, 1997, Lewis et al, 1999, 2000, Rios et al, 2013) and its generalizations (Kolda et al, 2006).

The most basic derivative-free local optimization method is coordinate search (or "ascent", or "descent"), in which each iteration directly tries to improve the results obtained by the previous iterate. The pattern followed in a search sequence is called a "sweep pattern". Coordinate search can be issued for one cycle of iterations; or can be repeated recursively without varying the evaluation order; or can be repeated while varying it: by inverting it at each cycle (Aitken method - Bazaraa et al, 1993) or by choosing more complex patterns of recombination. Search performance is known to increase when the set of search directions is rotated after each parameters have been used once (Rosenbrock, 1960); but there are smooth cases in which, by using coordinates as the only means of direction, convergence may fail (Powell, 1973). Convergence of coordinate search can be especially problematic for non-convex and discontinuous problems (Luo et al, 1992a, b), if countermeasures are not taken.

### 2.4. Optimization-aimed problem decomposition

A strategy to reduce the burden of global optimization is to subdivide global problems in problems structured in more manageable parts. Problem decomposition is today widely thought an essential part of problem modelling. It may be originated by subdividing a problem (1) horizontally, by resolution (Diao et al, 2011), or (2) vertically, by system (Jedrzejuk et al, 2002a, b, c) or component (Geyer, 2009); or (3) by expertise domain, like in Multidisciplinary Design Optimization (MDO) (Cramer et al, 1994).

*Dynamic Programming* (DP - Bellman, 1957) is a decomposition-based method especially suited to problems in which subproblems can be optimized independently. DP introduced the concept of memoization, i.e. the recording in memory of the results of recurrent calculations for subsequent use.

A stream of studies about decomposition has been motivated by distributed computing and the need of parallelization (Michelena et al, 1997, 1999).

The centrality of the issue of MDO-related decomposition for building design and the current lack of integrated CAD-centric tools allowing MDO has been stressed by Welle (2012).

#### 2.4.1. Block search methods

When the pursuit toward optima is operated by varying one subspace at a time, a search may be defined "block search", a descent "block descent", and a subspace a "block". The simplest formulation of block coordinate descent (BCD) is cyclic and based on fixed non-overlapping blocks. Search strategies internal to subspaces can be local or global and subspace searches can be linked to each other through local or global search strategies.

Block coordinate descent (BCD) is also known as block nonlinear Gauss-Seidel method (Bertsekas et al, 1989; Bazaraa, 1992; Berstekas, 1999). When there are overlaps between blocks, it resembles the multiplicative Schwarz alternating method (Schwarz, 1870; Saad, 1996). In literature, in the case of the Schwarz alternating method the term "subspace" is ordinarily used with respect to a partition unit of a problem, whilst in the case of the nonlinear Gauss-Seidel method, and less exclusively for BCD, the term "block" is used.

Cyclic block coordinate descent can be defined in the following manner: *min f(x)* subject to $x \in X$ such that $X = X_1 \times X_2 \times \cdots \times X_m \subseteq \mathcal{R}^{n_i}$. $X_i$ is a closed convex subset of $\mathcal{R}^{n_i}$ and $n = n_1 + n_2 + \cdots + n_m$. The vector x is partitioned in vector components $x = (x_1, x_2, \cdots, x_m)$, so that $x_i \in X_i$, for i = 1, 2, $\cdots$, m. The search algorithm solves for $x_i$ fixing and other subvectors of $X$ cyclically (35). The block which is tried at a certain iteration is indeed called the active (or working) block. Given an initial point $x^{(0)} \in X$, for k = 1, 2, $\cdots$, if $x^{(k)} = (x_1^{(k)}, x_2^{(k)}, \ldots, x_m^{(k)})$ is the current iterate, the next iterate which is generated is:

$x^{(k+1)} \leftarrow (x_1^{(k+1)}, x_2^{(k+1)}, \ldots, x_m^{(k+1)})$, according to the iteration:

$$x_i^{(k+1)} \in \underset{y \in X_i}{arg\,min} f\left(x_1^{(k+1)}, \ldots, x_{i-1}^{(k+1)}, y, x_{i+1}^{(k+1)} \ldots, x_m^{(k)}\right) \text{ (Berstekas 1989).} \qquad (1)$$

A formalisms to integrate in the descent algorithm the capability of allowing for overlapping blocks is presented in Grippo et al. (2011). In this formalism, at the kth iteration vector $\mathbf{x}^k$ is partitioned in two blocks, where $\mathbf{A}^{(k)}$ identifies the active set, composed by the variables that are updated, such that $\mathbf{A}^{(k)} \subset \{1, \ldots, n\}$, and $\mathbf{I}^{(k)}$ identifies the inactive set, composed by the other variables, such that: $\mathbf{I}^{(k)} = \{1, \ldots, n\} \setminus \mathbf{A}^{(k)}$. On the basis of the solution which is current at the kth iteration $x^{(k)} = (x_{A^{(k)}}^{(k)}, x_{I^{(k)}}^{(k)})$, the active block is calculated by solving the subproblem $\min_{x_{A^{(k)}}} f(x_{A^{(k)}}^{(k)}, x_{AI^{(k)}}^{(k)})$. The inactive block is not modified: $x_{I^{(k)}}^{(k+1)} = x_{I^{(k)}}^{(k)}$, and the current solution is updated by setting $x^{(k+1)} = (x_{A^{(k)}}^{(x+1)}, x_{I^{(k)}}^{(k+1)})$. \qquad (2)

The convergence of block-coordinate descent is problematic for non-convex problems and discontinuous problems (Razaviyayn et al, 2013). But several issues regarding BCD has been proven circumventable. Its convergence has been proven for both quasi-convex and pseudo-convex objective functions (Tseng, 2001; Razaviyayn et al, 2013). Moreover, the convergence of BCD in some kinds of non-differentiable problems has been proven when they can be reduced to dual problems (Tseng et al, 2009a, b).

In the case in which distinct accumulation points are present, a two-blocks method has been proven useful for convergence, obtained by keeping the search criteria for the second block less strict than those used for the first (Grippo et al, 1999, 2000). The application of that method to sparse problems through adaptive redefinition of blocks has also been studied (Grippo et al, 2011).

The quasi-equivalence of convergence criteria for non-overlapping and overlapping block coordinate searches has been proven (Chen, 2005; Cassioli et al, 2013; Razaviyayn et al, 2013) and the better convergence of random block coordinate search over cyclic block coordinate search for smooth and unsmooth problems has been observed (Nesterov, 2012; Richtárik et al, 2014).

**2.4.2. Considerations on problem decomposition in the context of building design**

Problem decomposition has been used in design well before computer existed. For instance, in cascading drawing scales in building design endeavours. Problem decomposition also occurs whenever a building simulation inquirer substitutes a thermal zone with another one modelled at a higher resolution. Some information in the two subspaces is repeated and is destined to stay unchanged; some information is repeated but is due to change; and some information is not shared. This is a very common situation, having the effect of lowering the probabilities that functional inefficiencies arise in designed and self-designed information-based systems, which include living organisms. That strategy has at its core the embodiment of a certain level of commonality, of shared information between instances (spaces, subspaces, or individuals).

Block nonlinear Gauss-Seidel, which has been adopted in the present research, is the simplest decomposition-based optimization method. Optimization and decomposition are tightly knit in it. It is a very general idea, which can be performed both by a man and a machine. The choice of that method leaves the emphasis of this study on decomposition.

**2.5. Metaheuristics**

The increase of computing power which has taken place in the last decades has prompted a rise of interest in global optimization methods. Those are ordinarily more useful than local ones for real-life design endeavours, which are often characterized by discontinuities and non-convexities hindering convergence towards absolute optima. And that interest has accompanied the growth of the so-called metaheuristics - i.e. low-level heuristics applied to search to make it more efficient (Blum et al, 2003; Glover et al, 2003; Gendreau et al, 2005).

The relative performance of some deterministic and probabilistic algorithms has been compared in the framework of building optimization by Wetter et al.[2004].

Amongst the most widespread metaheuristics are evolutionary algorithms (EAs), simulated annealing (Kirkpatrick et al, 1983; Černý, 1985) and swarm-intelligence methods (Kennedy et al, 1995, 2001; Colorni et al, 1991). The most widely used in are EAs.

**2.5.1. Evolutionary algorithms**

Evolutionary algorithms (EAs) are population-based and adopt the principles of natural evolution. Like direct search methods, they use discrete data, requiring that continuous problems are transformed into combinatorial ones. The first experiments with EAs were conducted by Nils Barricelli to investigate the role of symbiogenesis in evolution (Barricelli, 1954, 1962 a, b). Amongst the principal EAs are evolutionary programming (Fogel et al.,1966), based on finite-state machines; genetic algorithms (GAs – Holland, 1975-1992), closely inspired to the mechanisms of reproduction in living organisms; and genetic programming (Koza, 1992), in which the objects of evolution are the programs managing information themselves (so that evolution can also involve its own rules), organized in hierarchies of trees or graphs.

The most used EAs in design today are GAs, which adopt the full range of genetic operators seen in nature: recombination (deriving from crossover and inversion) and mutation. The choromosomes carrying the genetic information can be encoded (usually in binary format) or can be directly constituted by the design parameters.

The seminal interpretation of GAs given by Holland (1975-1992) stresses the importance of advantageous groups of genes (*schemata*) in chromosomes working as building blocks and constituting the minimum units of recombination in evolution. Advantageous schemata increase exponentially (in an infinite population) at each generation; to which corresponds the interesting property that disadvantageous ones disappear logarithmically, relegated into a sort of quiescence, of "reservoir" of unoptimal solutions which may result useful for giving a population the possibility of "backtracking" if it gets stuck into some local optima. That reservoir makes possible that the reproduction of individuals takes place like in a Markov chain (i.e. with no individual memory) without causing a global memory loss. The resilience allowed by this is further increased when (as made possible by meiosis and sexual reproduction) a schema is given the possibility of hiding itself in inactive (recessive) state in a genotype without harming an individual.

The necessity of "backtracking" may arise both when the modelled problems are static or can be assimilated into ones - like design endeavours performed once, which is usual for buildings - and when the environment keeps changing, so that (a) the past usefulness of a solution is not a guarantee for the future and (b) the design endeavour is "perpetual", like in biological life. The consequence of this "dilution" of memory in gene pools is that GAs are characterized by the need of a balance between exploration (driven by recombination and mutation) and exploitation (driven by selection).

Many GA variants have been experimented; among them, elitist GAs, in which the ordinary genetic operators are accompanied by additional selection criteria based on performance; which makes them most competitive for problems also suited to enumerative optimization methods - to which they may be considered in a sense hybridized. The most widespread elitist GA is currently the NSGA-II (Deb et al, 2002).

## 2.6. Relations between evolutionary algorithms and decomposition-based enumerative methods

The concepts of EAs and block coordinate search (BCS) have areas of overlap. BCS may be viewed as a particular case of EA in which: (a) the population is composed by one individual; (b) genes are not encoded; (c) mutation is (usually) determined on the basis of enumerative criteria; (d) the "parents" are constituted by the "survivor" of each block search phase (i.e. generation) and the individual existing before that phase; (e) the crossover area corresponds to the active block; (f) when the active block is "donated" in a crossover-like operation, the homologous "chromosome" is not formed. Moreover, (g) blocks in BCS can be defined so that they correspond to the groups of *loci* defining schemata in EAs, with the consequence that it is possible to design genetic algorithms having a decomposition structure devised by expert knowledge.

What ultimately distinguishes BCS from EAs is that the former lacks crossover between coexisting individuals. In BCS different individuals can exist in the same time and even pursue different objectives, but they cannot exchange information, and cannot therefore profit from the information embodied in coexisting individuals. Which implies that searches lack a "storage unit" of information about their own history, a "latent memory". This marks an advantage of GAs over BCS.

The second condition deriving from not being population-based is that in BCS the current individual is always the fittest ever recorded, so that there is no possibility of "backtracking" from local optima by returning on previously on-the-way-to-be-discarded solutions. The current individual can only escape from local optima by sidetracking towards more nearly optima. An interesting consequence of this is that recombination in BCS does not implies disruption; which avoids the tradeoff between exploration and exploitation typical of EAs and allows for more daring and "enterprising" explorations. Which is an advantage of BCS over EAs.

**2.7. Search in building design**

Simulation-assisted optimization strategies in building design are today mostly used for fine-tuning in specific, later-stage types of inquiry (Augenbroe, 2003), characterized by relatively clear goals and relatively few variables in play; or targeted to specific, constrained types of problems, determined by building type, or use, or objective function, and by the selection of specific issues of inquiry (Daum et al, 2009; Shin, 2011; Khan et al, 2012). This happens in spite of the fact that the variables influencing building performance are numerous and interrelated. Parametric analyses today indeed seldom regard morphology, or involve only more or less narrow subsets of morphological variations (Caldas et al, 2002, Wright et al, 2002; Djuric et al, 2007; Gong et al, 2012; Ihm et al, 2012). And issues of shape are usually confronted through very specialized, dedicated models and schemes (Adamski, 2007; Turrin et al, 2011). But experience shows that morphology is the all-dominating factor in architectural design, also due to its influence on space usage and perception.

The fact that the least assisted design phases are the early ones contributes to push architectural designers toward a tentative kind of approach during the early phases of design, often based on almost blind samplings in problem space. As a result, contemporary projects often carry the marks of premature optimization, due to a hasty exploration of options early on combined with a higher focus at later stages. It is likely that historically the risk of premature optimization in building design got worse with the possibility - and the need - to choose; specifically, with the design and building power brought by modern methods of production and the scientific method, long evolved through history (White, 1964). Building had grown more sophisticated after that, but taking decisions had become increasingly difficult.

Symptoms of premature optimization were indeed rare in traditional architecture, whose solutions (although clearly much more constrained) were selected by time and use whilst remaining generic enough to be robust, applicable in a variety of situations (Rapoport, 1969). That shift has also been due to the increase in construction speed brought by modernity. Pre-modern construction techniques were slow. As a result, a great deal of design decisions were differed to the on-site, construction phase (Viollet-le-Duc, 1873). This let designers dwell on problems (indeed, it required they did), which was beneficial for the quality of design (Alexander, 1979), often to the point of making solution context-aware in spite of a lack of scientific knowledge. A "design as you build" approach may also be adopted for buildings constructed with modern methods, but it is likely to produce abrupt changes in the process (Brand, 1994), to the point of stimulating researches on modern technologies suitable to allow on-site decisions (Alexander et al, 1985).

A major cause for today's scarcity of approaches supporting the early stages of design is likely to be that the capability to morph models along flexibly specifiable paths of modification is difficult to obtain from the combination constituted by parametric exploration tools and building performance simulation (BPS) tools. An important reason for this is that the capability to morph models along flexibly specifiable paths of modification is difficult to obtain from the combination constituted by parametric exploration tools and building simulation tools: today's state-of-the-art tools for parametric analysis, both building-specific (Mourshed et al. 2003, Christensen 2006, Caldas 2008, Zhang 2012, Palonen et al. 2013, Ellis et al. 2006, Attia et al. 2012, 2013) and multi-purpose (Adams et al. 2011, Wetter 2000a, b), require a more or less explicit description of all the "actions" to be performed on models, which may be long and difficult if the actions are complex and intertwined.

Constraint-based morphing procedures can raise the level of abstraction and flexibility of the description of model mutations because they treat them relationally, taking into account how all the parts of a model "react" when some of them are changed. The approach based on propagation of constraints is present since decades in both research (Brüderlin, 1986) and the programming scenery (Chenney, 1994) and has been integrated into simulation-aimed models (Yi et al, 2009). The technical barriers impeding the objective of bringing dynamic properties to models are being removed and that objective is being pursued through parametric scripting (Nembrini et al, 2014), integration with BIM software (Gerber et al, 2013; Welle et al, 2011) and the recourse to shape grammars (Granadeiro et al, 2013).

The choice of algorithms for building optimization is determined by the fact that the problems to be confronted are mainly non-linear and multi-objective. It is therefore no surprise that recent studies (Kämpf et al, 2010 a), reviews (Evins 2013, Nguyen et al, 2014) and experiences (Yi et al, 2009; Magnier et al, 2010; Turrin et al, 2011; Evins et al, 2012) show a stronger and stronger prevalence of metaheuristics; mainly, GAs - with a prevalence of the elitist NSGA-II.

**2.8. Search in micro-urban design**

At micro-urban level, building and micro-urban simulation tools converge.

When building-scale tools are used, the problem can take the form of the search of the position in a site, requiring that the candidate positions are sampled one at a time, due to the limited space that can be covered, and the parameters describing position are assumed as variables in problem space, with the result of increasing its size and adding to it some degrees of randomness and discontinuity deriving from the urban context. In those situations, the capacity of exploring the search space globally is most needed when the elements of randomness and discontinuity typical of real architectural scenes - generated by buildings, trees and shrubs, terrain shapes, water bodies - are brought into description. In building simulation tools these entities may be represented as solar and wind obstructions, solar reflectors, and far-infrared emitters (Crawley et al, 2008); but only the first two are taken into account when open spaces are modelled as external to the assessed thermal zones. This can lead to underestimations of discontinuities and non-convexities in the objective functions and can be confronted by modelling the external spaces themselves as thermal zones. That strategy produces a better predictive capacity, but at the cost of a greater modelling complexity; and is in need of more research about the criteria for modelling the outer, fictitious boundaries.

For what concerns the micro-urban scale, solar irradiation can be studied with the same tools suitable for the building scale (Compagnon, 2000), and simulation tools specifically conceived for the micro-urban scale are consolidating, focusing on solar irradiation (Dogan et al, 2012), building energy consumption (Robinson et al, 2009), and thermal radiant exchanges (Hénon et al, 2011). CFD-based tools (Huttner et al, 2009), are becoming more attractive for parametric studies due to the reduction of their computation times; and a new class of tools integrating all the above kinds of analyses is beginning to appear (Musy et al, 2014, Gracik et al, 2015).

As confirmed by recent researches (Oliveira Panão et al, 2008; Kämpf et al, 2010a, b; Martins et al, 2014; Okeil, 2010) and reviews (Srebric et al, 2015), micro-urban-simulation-based searches make possible a simultaneous analysis of all positions in space, which avoids the necessity of including the parameters describing position as variables in problem space. This keeps the problem size smaller, but doesn't obviate to the elements of discontinuity and randomness which are present in those situations. For this reason, those searches as well are currently dominated by GAs.

### 3. Definition of a parametric test-case for the experiments

A building-simulation-assisted test-case was set up to study the implications of different search structures for micro-urban design. The optimization process regarded the features of the building (construction and form) and its position in space and was pursued on the basis of environmental performances.

In a certain climate (Milan) a space was defined by placing solar and wind obstructions. There, a small single-zone building with a given volume had to be situated and designed. The design parameters were 10, regarding the building position (x and y coordinates), its shape (about which width, depth, convexity, rotation; and front, side and back windows surface percentage were to be defined) and the building envelope's effective thermal capacity (Fig. 1). The chosen parameters were tested for 3 values each, producing a total search size of $3^{10} = 59049$ evaluations.

After the given changes to the model were executed, insolation, shading and flow network information were updated to take into account the modifications of building and urban obstructions regarding solar radiation and wind flows. (Further information in the Appendix.)

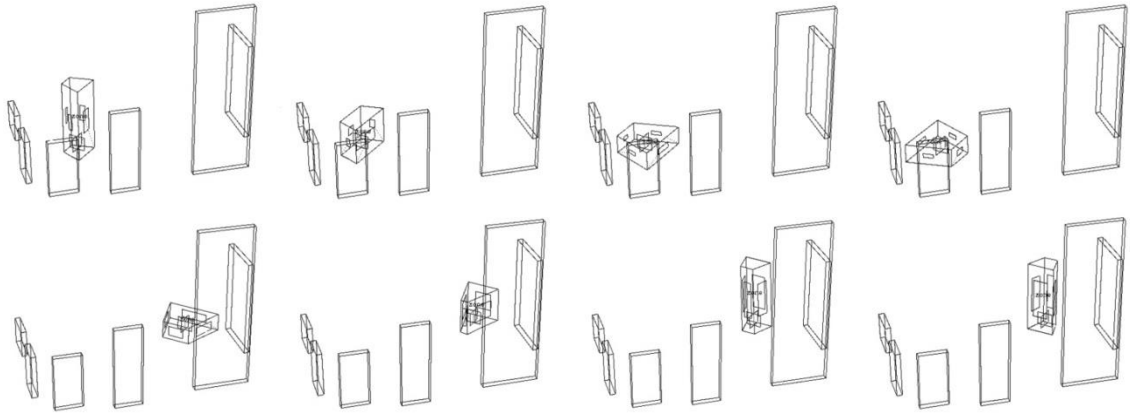

Fig. 1. Some sample instances in the morphing sequence. The shape of the building is the one mutating at the centre of each scene. The boxes around it are solar and wind obstructions.

Two objective functions were taken into account: winter (February) heating loads and summer (August) average maximum resultant free-floating temperatures (chosen to verify overheating avoidance). They were normalized by scaling them from 0 to 1; then a similarly normalized bi-objective functions was derived by weighting them 1:1. The mono-objective functions were chosen to be partially conflicting. This is because site obstructions reduced winter and summer solar gains with a different timing depending on the season. Being the site at a median latitude (45°), obstructions at the south edge of the plot were more influential in winter and obstructions at the east and west edges in summer. Conflict could arise seeking for a position because the sizes of the two side windows varied together. When the building was placed near the west or east of an obstruction, a conflict could arise in summer, because the window near the obstruction was more sun- and wind-shaded than the opposite one. As a result, the distribution of performances deriving from variations in side windows was non-convex (Fig. 1, rightmost graph, upper row).

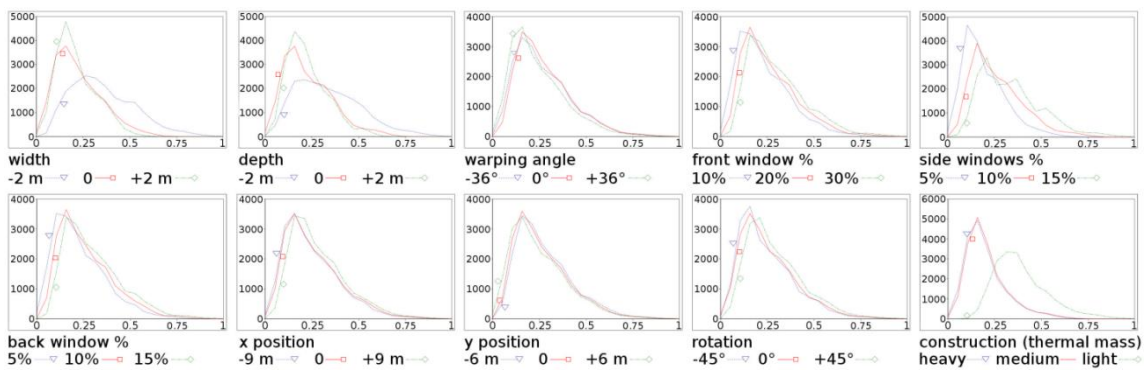

Fig. 2. Distribution of performances generated by a linear variation (brute force optimization) of the 10 parameters assuming 3 values each. On the abscissae: value of the multi-objective function. On the ordinates: number of cases.

### 3.1. Implementation choices

The mutation of models were managed through the use of the *Sim::OPT* morphing program (Brunetti, 2008-2015), which pursues relational capabilities by propagation of constraints (Abelson et al, 1985-96) targeting the *ESP-r* building performance simulation platform (Clarke, 2001; ESRU, 2001). Sim::*OPT* works by sending inputs to the *ESP-r* user interface (Hand, 1998). After each mutation event, it reads the *ESP-r* model description files; gets information back by reading the *ESP-r* model's databases (ESRU, 2002; Citherlet et al, 2002); checks if the model satisfies the user-specified constraints; and, if needed, transmits requests for modification of model files back to the BPS platform. Propagation of constraints eases the separation of the description of the events triggering morphing from the morphing rules (Brunetti, 2013).

### 3.2. Definition of commonality

The global search strategy that was used in subspaces is the most basic one: linear search. The active blocks are searched one after the other, one at each morphing step, while the inactive block is kept fixed; and the most performant parameter values (vector components) are selected and assumed for the following iterations, so that the objective function is always minimized or at least unvaried. Due to the overlaps, the parameters belonging to the intersection set between subsequent active blocks do not pass values: they are allowed to vary so that they cannot exert a selective effect over the values of the parameters to be passed.

When active blocks (subspaces) overlap and the evaluations which are shared between them are repeated, duplicated information occurs, and redundancy with it. But if the duplication of evaluation is avoided, thanks to memoization, like it was done in the present study, redundancy can be viewed as commonality; that is, shared structure without repetition. The term "commonality" will be here used from now on to indicate shared structure (overlap) purged from redundancy through memoization.

## 3.3. Choices regarding search structures

Sweep patterns were made dependent only from the iteration numbers and search cycles. The effects of sweep patterns which were different as regards both type and duration were evaluated. In Figures 3, 4 (and 5) each box represents one parameter and each arrow represents the information flow determined by the passing of minimizing coordinate values from one iteration (k) to the next (k+1).

In the experiments, all searches were tested on the basis of the same 590 random sequences of parameters, and the considered results were always the average results of these 590 cases. This was done to isolate the influence of the search structures from that of the assignments of parameters to subspaces. 32 mono-cycle searches (7 shown in Fig. 2) and 54 multi-cycle searches (Fig. 4, 6) were tested. The 54 multi-cycle searches were derived from 18 decomposition schemes (Fig. 5) and 3 sweep patterns (Fig. 4). Preliminary inquiries were performed to select an approach suited to the objectives and an adequate block search method. As regards the modelling of objectives, after evaluating the effects of penalty and parabolic functions, an ordinary linear weighting was chosen. As regards the search method, a Gauss-Seidel method (sequential) was preferred to a Jacobi one (parallel), for the sake of search quality.

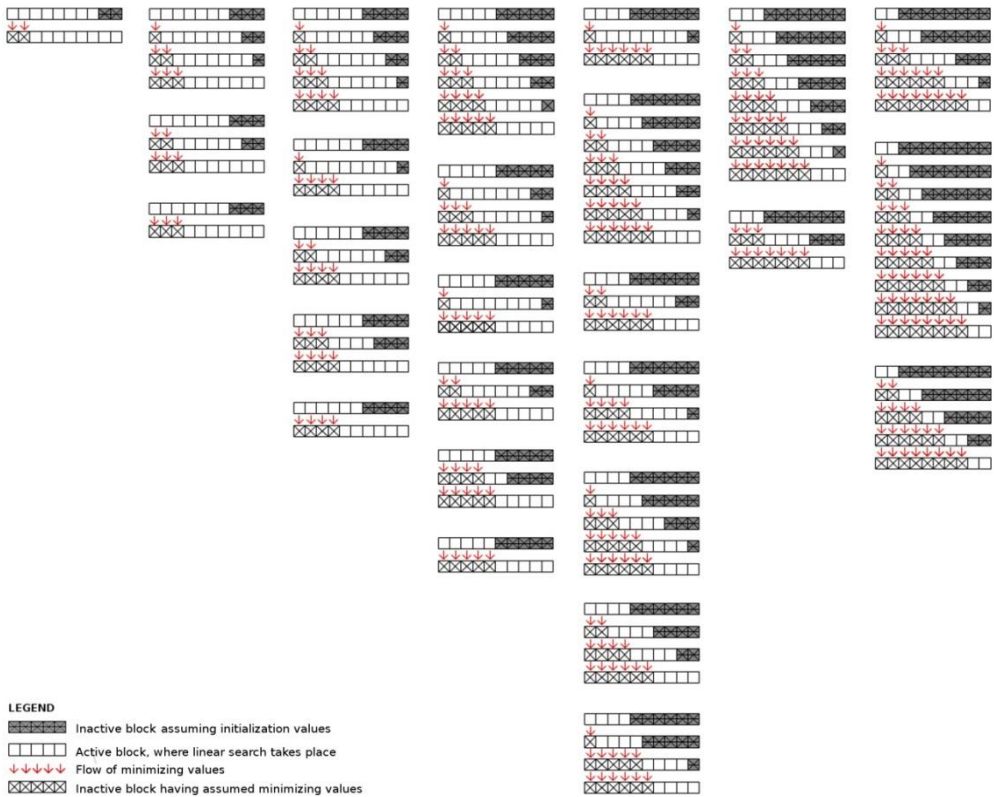

Fig. 3. Schemes of the search structures (overlapping block decompositions and sweep patterns) tested in the mono-cycle, non-recursive searches that are shown in Fig. 6.

**3.4. Contents of the main tests**

The following aspects were inquired:

1) the role of commonality between subspaces;

2) the effectiveness of search repetitions without sweep pattern variation;

3) the usefulness of recombination of subspaces coupled with repetition.

The search structures that were adopted to test the recombination types in Fig. 4 are shown in Fig. 5. In principle, successive sweep cycles would not necessarily need to follow regular patterns, so that even the concept of sweep cycle may result unfit; and active blocks would not need to be constituted by adjoining parameters. But those settings have been assumed to ease the isolation of conditions determining search quality and efficiency.

**3.5. Criteria for recursive recombination**

The recombination of subspaces required the possibility of traversing the lateral boundaries of search spaces (Fig. 4). This suggested the adoption of circular lists, which could make possible for a subspace exceeding the tail of a list to re-enter it from the head, and vice-versa. The required functionality was obtained by joining copies of Perl arrays and operating at the middle of the so-obtained sequences.

Together with search iteration (column on the left in Figure 4), two kinds of recombination were tested. Both obtained recombination through an inversion of space traversing verse at each sweep cycle and defined the starting element of each sweep cycle as the farthest non-visited mid-point. In recombination type "A", the list positions of the first varying element at each sweep cycle were 0, 9, 5, 4, 7, 6, 2, 1, 3, 8. In recombination type "B", they were 0, 5, 7, 2, 9, 4, 8, 3, 6, 1. Type "A" in most search structures produced a fuller recombination. In the tested searches, it never happened, however, that it took more than 6 search cycles before the improvement of search quality wore out (Fig. 9, 11).

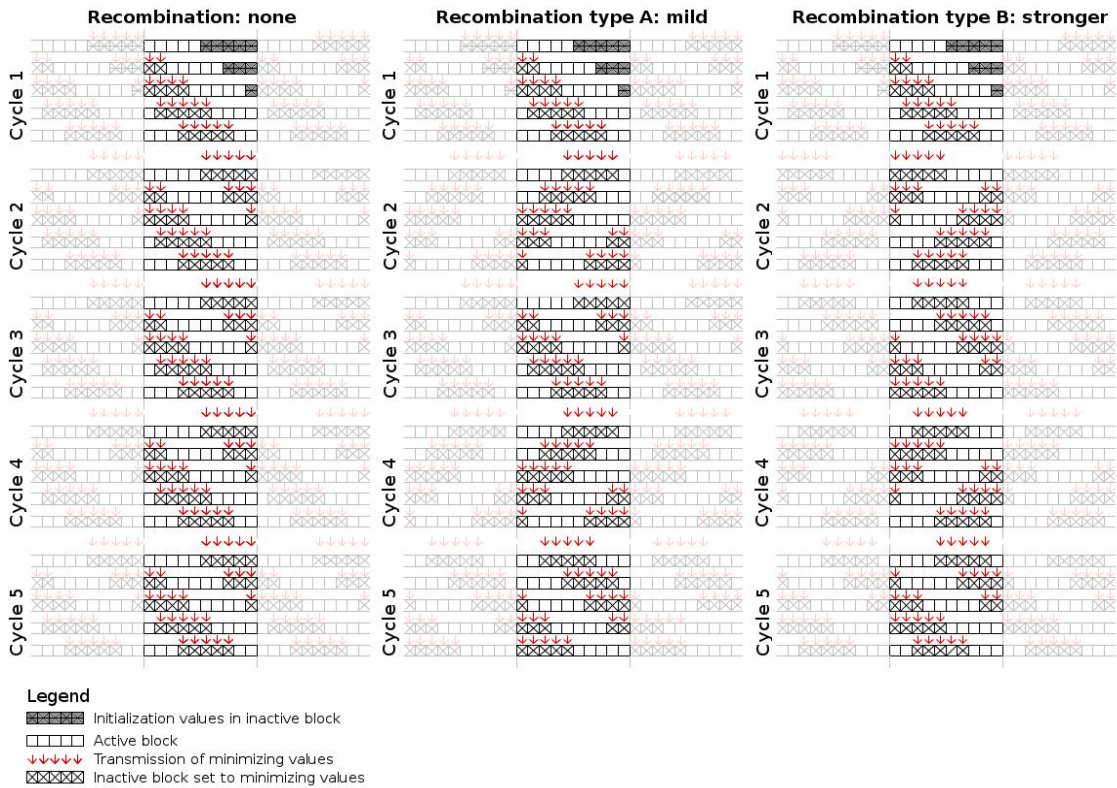

Fig. 4. Three example schemes of overlapping block coordinate descent featuring repetition plus recombination and modelled through circular lists. The considered types of recursive search patterns are: repetition without recombination, repetition plus recombination (A), repetition plus a stronger recombination (B). In strong black: the sequences (lists) with which the parameters in the search space are modelled. In pale grey: the blocks that enter or exit the lists from their head and tail.

| Active block size | 5 | 4 | 3 | 2 |
|---|---|---|---|---|
| **Continuous search structures** — Overlap size 4 | B5-O4 | | | |
| 3 | B5-O3 | | | |
| 2 | B5-O2 | B4-O2 | B3-O2 | |
| 1 | B5-O1 | B4-O1 | B3-O1 | B2-O1 |
| 0 | B5-O0 | | B3-O0 | B2-O0 |
| **Truncated search structures** — 4 | T-B5-O4 | | | |
| 3 | T-B5-O3 | | | |
| 2 | T-B5-O2 | T-B4-O2 | T-B3-O2 | |
| 1 | | | T-B4-O1 | T-B3-O1 |

Legend: ▭▭▭▭ Active block, ▨▨▨▨ Inactive block

Fig. 5. Recursively recombined overlapping block coordinate descent schemes tested in the trials. Only the first search cycle is shown here, for brevity. The recursive parts of these search structures follow the patterns shown in Figure 4. These sequences are circular as well.

### 3.6. Definition of heuristic indicators measuring commonality

The role of commonality was investigated by the means of indicators defined on purpose. The first one - *commonality ratio* (CR) - measured commonality as the ratio of the overlapping block search size and the total search size (Table 1). From commonality ratio (CR), *commonality volume* (CV) was obtained, as the product of CR and net search size (SS - Table 1). But CR had the limit of not being sensible to the bottlenecks that are present in search structures when an overlap between subspaces is absent or small. To deal with that situation, a *local commonality ratio* (LCR) was defined (Table 1), constituted by the ratio of the cumulative search size of two subsequent active blocks (subspaces) and the search size shared between them (commonality). From it, a function that will be here named *commonality flow* (CF) was derived, constituted by the cubic root (smoothing variations) of the cumulative product of the active subspace size and the local commonality ratio at each iteration step (Table 1). Being a cumulative product, this function is capable to signal the effect of bottlenecks in a search structure and can be used to measure its level of integration. From CF, as its cumulative sum, an indicator that will be here named the *cumulative commonality flow* (CCF) was obtained. Being a cumulative sum, it can only increase with iterations, as information builds up.

## 4. Results and discussion
### 4.1. Non-recursive searches

The role of commonality between subspaces resulted to be positively correlated to search quality in one-cycle, non-recursive searches (as shown by CF, Fig. 2, 6). A high commonality between subspaces appeared to produce a high search quality. This must be due to the information flow it allows. Commonality, produced by overlaps, warrants indeed a good "welding" between subspace searches. For this reason, in devising the considered search structures, it may be useful to place the most responsive parameters in the overlaps. However, a high commonality produces a great search size, which comes at the expense of search efficiency.

The behaviour of those indicators can be observed in Figure 6, where the mean results the 590 considered sequences relative to some mono-cycle, non-recursive searches are shown in order of decreasing search size. In the adopted naming criterion, B7-O6, for example, indicates a search structure composed by active blocks being 7 parameters in size and overlapping for 6 ("B" stands for active block and "O" for overlap). Had the sweep pattern been truncated, its name would have been preceded by a "T-".

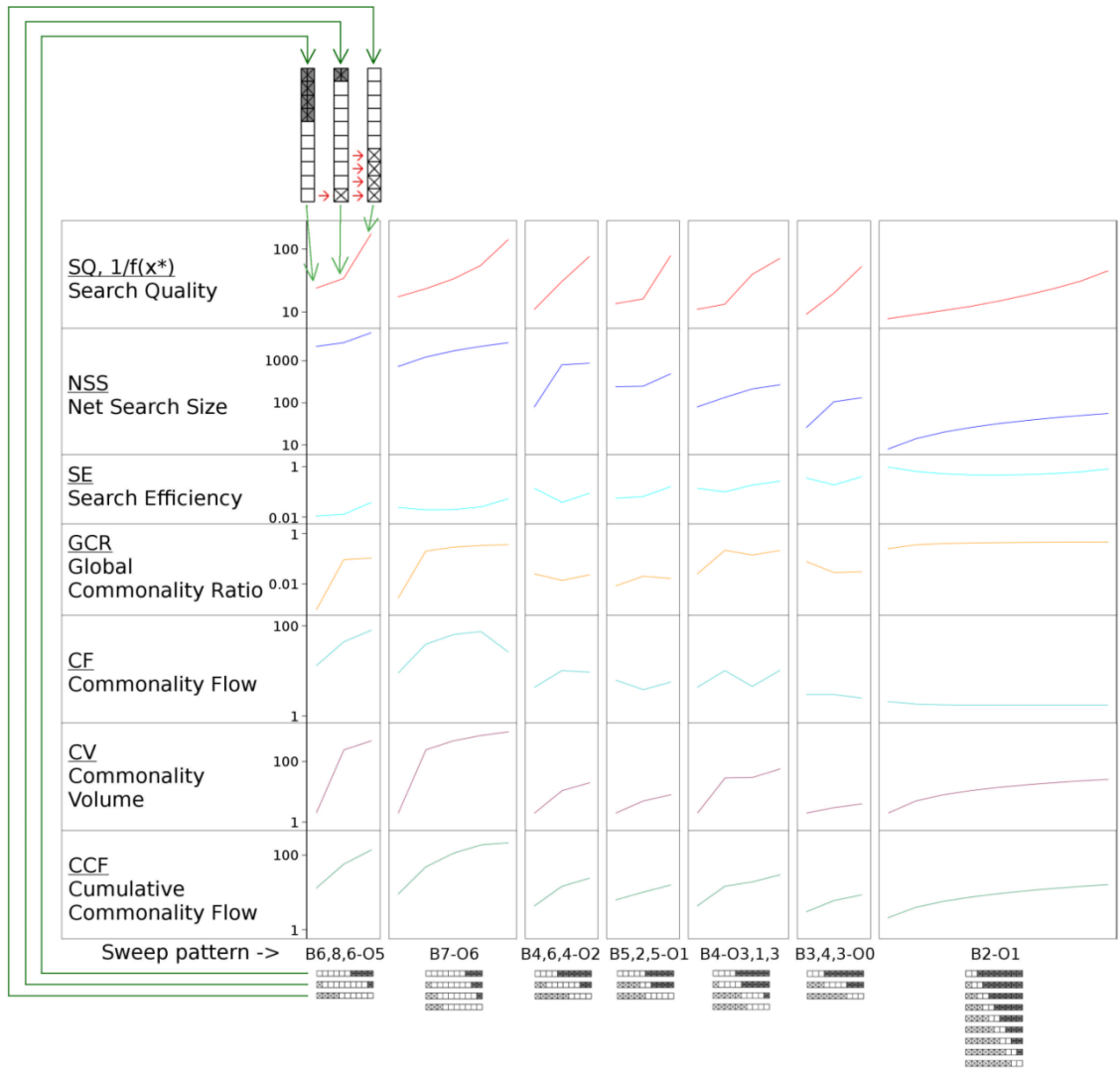

Fig. 6. Sample results of the mono-cycle searches shown in Figure 3. The minimization problem was transformed in a maximization one, for clarity. The scales are logarithmic. Above on the left, the links between subspace searches and plotted results for a search structure are highlighted.

Because the searches in Fig. 6 were started with initialization values, only the results obtained in their last iteration are intended to be meaningful as for performance. From that fact also derives that search efficiencies in the plots describe upward concave curves, recovering toward completion.

From these results it can be verified that whilst it is possible to compromise the search performance of a search structure of large size due to a poor internal integration, little can be done to improve the search performance of a well-integrated non-recursive search of small size. In the shown cases, search quality mainly vary with search size, global commonality ratio, commonality volume and cumulative commonality flow. And as a consequence, it never happens that a high search quality takes place together with a high search efficiency: they exclude each other. It seems difficult for a non-recursive, non-recombined search process to escape the fate written in its size.

Cases B6,8,6-O5 and B7-O6 are of big size and converge well; cases B4-O3,1,3, B3,4,3-O0 and B2-O1 are of small size and perform worse. When the size of such a non-recursive process is known, the most about its potentialities is known. For its good correlation with search quality and its dependence on search structure, CCF seems to be suited to support the design of non-recursive search structures; for instance, to choose among candidate structures of similar size.

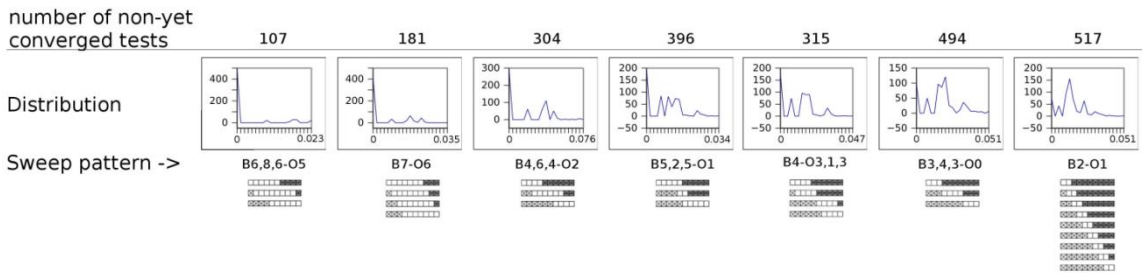

Fig. 7. Distribution of search performances regarding the non-recursive searches shown in Figures 3 and 6. The problem is shown as a minimization one (best results are on the left of the scale). On the ordinates: number of cases. On the abscissae: value of the multi-objective function.

Size being equal, searches characterized by higher CCFs showed indeed to produce higher search performances. This is what happens for instance in case B5,2,5-O1, which outperforms B4,6,4-O2, of greater size.

Highly integrated search structures, like B8,6,6-O5 and B7-O6, resulted to be more likely to produce optimal results. The plots show indeed that the more integrated a search structure is, the greater is the number of search sequences converging to optima. Possible causes of low integration may be small search sizes (cases B4-O3,1,3 and B2-O1) or bottlenecks in the search structures (cases B4,6,4-O2, B5,2,5-O1 and B4-O3,1,3).

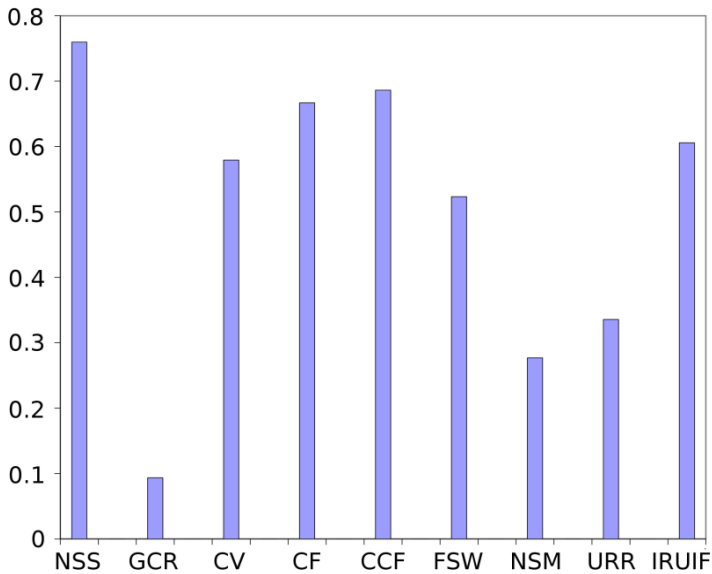

Fig. 8. Correlation between some of the functions regarding mono-cycle, non-recursive searches, normalized by scaling the maximum to 1, and search quality ($1/f(x^*)$). Legend: NSS: Net Search Size; GCR: Global Commonality Ratio; CV: Commonality Volume; CF: Commonality Flow; CCF: Cumulative Commonality Flow; IRUIF: Indicator of Recursive Usefulness of the Information Flow.

### 4.2. Recursive searches

#### 4.2.1. Search iteration

Search iteration demonstrated to be capable of improving search quality, the more so the less the first sweep cycle was integrated internally (i.e. the less overlapped its subspaces were), with decreasing results as repetitions went on, and a less rapid wearing out of improvement in convergence the smaller the search integration and the higher the search recombination (Fig. 9, 11).

Now, because search repetition contributes to propagate information along an incompletely exploited block coordinate search, it can be used in place of commonality for that purpose. To adopt an analogy with a fluid flow, a large search size coupled with a high commonality may be seen as a large pipe, which allows a great flow even if the velocity is low. Search repetition is instead like allowing more time for a fluid to pass. But repetition when used alone as a self-sufficient strategy is exposed to a higher risk that a search does not converge, getting trapped into local optima. This limits the usefulness of repetition, if repetition is not accompanied by devices to counter the entrapment.

An examination of the ratio of search quality to search size, which defines search efficiency (Fig. 9), shows that repetition is useful to improve it, but especially when commonality and search size are low (B5-O0 A and B, B3-O0). In those cases, repetition produces a more considerable improvement of search quality, and it takes a greater number of repetitions before the improvement wears out.

### 4.2.2. Recursive recombination

Results showed that to improve search quality through repetition, recombination of subspaces at each search steps is essential, at a degree depending on the type of recombination. The most advantageous of the three tested recursive search strategies resulted to be the one adopting alternated sweep directions and a higher recombination of active subspaces; and each type of recursive recombination showed a peculiar amount of capability to improve search quality and efficiency: the fuller the recombination, the higher the search efficiency. Recombination coupled with search repetition resulted to be an alternative to commonality and search size to increase search quality in multi-cycle, recursively recombined searches, because it made less likely that searches did not converge by getting trapped into local optima. Repetition coupled with recombination produced higher search efficiencies and search quality (see cases see cases B5-O0 and B3.O0 in Fig. 9, and cases T-B5-O3, B5-O0, B4-O1, T-B4-O1 and B3.O0 in Fig. 11).

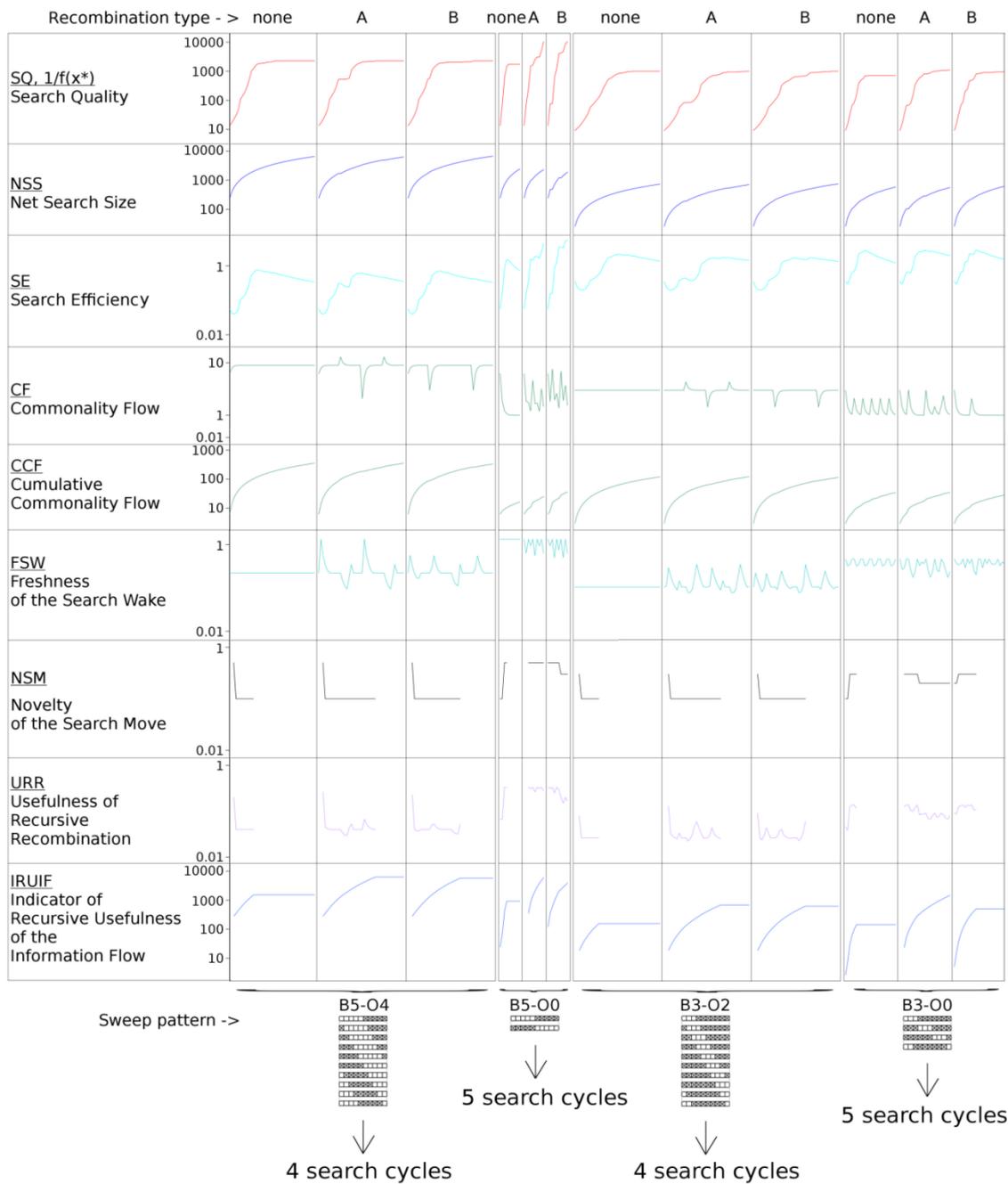

Fig. 9. Sample results of the searches shown in Figures 4 and 5 regarding the multi-cycle, recursive cases. The minimization problem has been transformed in a maximization one. The scales are logarithmic.

The importance of recombination to avoid that a search process is trapped into local optima is evident in Figures 9 and 11. Figure 10 shows a substantial improvement in the distribution of results obtained by the same search structures dealt with in Figure 9 at the final sweep iteration (Fig. 7). It can there be noted that the lack of convergence of certain combinations of parameters is due to the fact that the searches get trapped in some local optima. Recombination had untrapped most of them in the preceding sweep cycles; but new, untried moves would now be necessary to unlock the still locked searches. The smallest likelihood of entrapment happens in the cases in which a search is not overexploited before it is recombined, like for instance in the well-performing B5-O0. Overexploitation (due to too tight overlaps and a too homogeneous search structure) took place instead in cases formed by large and highly overlapped subspaces, like B5-O4, where searches ended trapped early.

Too frequent entrapments into local optima are resulted to be the main cause of poor mean performance in all the tested cases. Being the causes of entrapments often not known or confidently foreseeable in ordinary optimization problems, it is likely that a sound criterion for the design of such searches is that of choosing search structures capable of reducing the risk of entrapment in the first place, irrespective of parameter sequences (i.e., irrespective of the assignment of parameters to search structures).

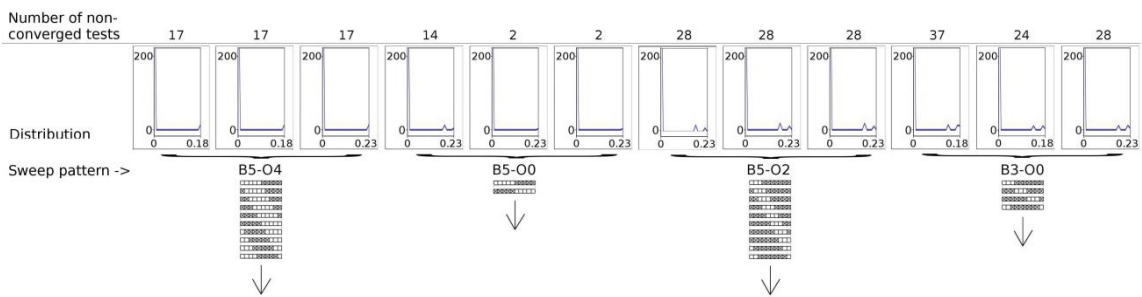

Fig. 10. Distribution of search performances in the last sweep iteration of the recursive structures shown in Figure 9, tested for the 590 random sequences. The problem is here shown as a minimization one and the converged tests are on the left of the scales. On the ordinates: number of cases. On the abscissae: value of the multi-objective function.

Recombination after a first search cycle appears to produce a real commonality between subspaces. But it also adds irregularity to the information flow, which warrants a more thorough exploration of subspaces, for the sake of convergence. Recursive recombination may be seen as originating a dynamically changing commonality. As an effect, repetition coupled with recombination postpones the wearing out of search quality improvement. And since it does not change the quantity of information flow, this must be due to the fact that it increases the quality of information flow, in the sense of its usefulness for search.

Recombination coupled with repetition can be of help to reduce the growth of search size deriving from an increase of the number of parameters, because by reducing the need of overlap between subspaces, it allows to reduce the size of subspaces. This is evident from the fact that search efficiency is high in searches characterized by a low internal integration and intense recombination, like in T-B5-O2, B5-O0, T-B4-O2 in Figure 11.

Recursive recombination makes also possible that the truncation of a search structure results in improvements in search efficiency. This happened for the truncated structures that were tested (T-B5-O3, T-B4-O2, T-B3-O2, T-B3-O1, Fig. 9), which is likely to be because whilst truncation reduces search size and commonality, it allows for the dismissal of search paths before they are fully exploited. In those cases, the obtained search efficiency results to be reduced at the first sweep cycle, but increased after a few repetitions.

### 4.2.3. Definition of indicators measuring recursive recombination

On the bases of the described observations, indicators were defined linking search size, commonality, number of search repetitions, number of recombination iterations, type of recombination, search quality and search efficiency. These definitions turned out to be necessary since the cumulative commonality flow (CCF) did not show a high correlation with search quality and efficiency in multi-cycle, recursive searches (Fig. 9). The key issue in the described framework is the definition of adequate indicators for recombination. The standard deviation of commonality flow was discarded for that because it would have implied that recombination could only take place in search structures characterized by alternations of bottleneck and enlargements. To quantify recombination, a more precise idea of the trace of past searches has to be sought. The concept of *search wake* has been here defined for that purpose; the search wake being constituted by the subspaces visited in the near past of a search. The search wake can therefore be monitored to gather information about the thoroughness of a search action. A useful kind of information can be obtained from the average age of the most recent parameter instances in the search wake, measured in terms of subspace (i.e. iteration) cycles. The ages of the cycles has been calculated by assigning a weight of 1 to the current cycle and kth to the first.

The obtained *average age of the search wake* (AASW) is a measure of how "fresh" a search wake is. Which is useful in the assumption that the more intense a recombination process is, the fresher the search wake is. The reciprocal of the AASW, which will be here termed *freshness of the search wake* (FSW), showed the highest correlation with search quality and efficiency for non-recursive searches (Fig. 8) and a good one for recursive searches (Fig. 12). Indicators alternative to FSW were tested, like the *refresh*

*ratio of the search wake* (depending on the length of a sweep cycle and AASW); but FSW was selected as more general, in that independent from the presence of regular sweep cycles, and even from the concept of cycles other than that of subspace search cycle; which makes that indicator also suited to deal with random searches.

The freshness of the search wake (FSW) is a measure of how quickly and thoroughly a search structure is traversed by subspace searches, but it cannot predict any consequences arising from recombination that do not alter the age of the search wake. This can be verified in the less recombined cases (the first ones in the plotted order) relative to search structures B5-O0 and B3-O0 in Figure 9, where it can be seen that the most recombined search structures (subspace size and sweep pattern duration being equal) obtain the highest performance. But FSW cannot explain that. One more indicator is therefore needed to predict the consequences of recombination.

That indicator, which will be here named *novelty of the search move* (NSM), is based on the information describing what ratio of each move from an active subspace to the next one is novel with respect to the past moves. This information is useful in the assumption that not only knowing if certain parameters were already visited is important; rating at what extent the path moves leading to the current active subspace are novel, never made in the search history, is important as well.

The ratio constituting the novelty of the search move (NSM) is calculated by describing the search moves from one subspaces to the next - both modelled as numerically ordered lists - as a sequence of pairs, each formed by an element of the previous active subspace and an element of the current one in the same position (following a certain order, in this case from the leftmost to the rightmost element, taken in ascending numerical sequence; then by counting how many of those pairs were chanced compared to each of the previous search moves. The highest value is selected; and finally the so obtained value is divided by the total number of parameters taken into account.

The FSW and NSM can be used together while keeping them separate or can be compounded into one indicator. That indicator will be here named *usefulness of recursive recombination* (URR), and calculated as the square root of their product.

Cumulative commonality flow (CCF) and usefulness of recursive recombination (URR) - or freshness of the search wake (FSW) and novelty of the search move (NSM) in place of the latter - are to be used together for the design of recursive search structures. This is because the good working of a recursively recombined search depends on a balance between commonality and recombination, which can be in conflict with one another. Commonality is indeed likely to be high when the overlaps between subspaces are large; but recombination may be high when the overlaps between the present subspace and the near past ones are small. To compound the above indicators into one, an *indicator of recursive usefulness of the information flow* (IRUIF) was defined as the cumulative sum of the product of CCF and the square of URR.

**4.2.4. Indicators measuring recombination**

The results in Figures 9 and 11 are presented in order of decreasing search size. In Figure 9, the performances produced by four search structures are reported, chosen in that representative of the general results: B5-O4 (high size, low CRF); B5-O0 (mid-high size, high CRF); B3-O2 (low size, low CRF); B3-O0 (low size, high CRF). Each of the four cases is analysed with respect to the three tested recombination types (following Figure 6: absent; A: mild; B: stronger, as signalled at the top of the figure).

The indicator of recursive usefulness of the information flow (IRUIF) showed the highest correlation with search quality and efficiency for recursive searches among the considered indicators (Fig. 12). It can then be observed that CCFs and FSWs are likely to become irregular when recombination is present. Then it can be noted that high FSWs, NSMs, URRs and IRUIFs appear to be the condition for the obtainment of high recursive search qualities, like in cases B5-O0-A, B5-O0-B and B3-O0-A and B3-O0-B. These searches are intensely recombined, and in them a high search quality takes place together with a high search efficiency.

Searches characterized by low FSWs and NSMs, like B5-O4 and B3-O2, reach lower search qualities and lower search efficiencies. And in cases characterized by low FSWs, NSMs and IRUIFs, search efficiencies top earlier than in searches characterized by high FSWs, NSMs and IRUIFs, which improvement lasts longer and goes further.

The correlations that was ascertained suggest that the performance of recursively recombined block searches depends on both their internal integration (evaluated through CCF), amount of subspace refreshments (evaluated through FSW) and amount of subspace recombination (evaluated through NSM), and that those conditions can be heuristically measured. Recursive searches showed indeed to behave differently from non-recursive ones. In particular, they showed that if a search process is characterized by a high recombination, it can escape the fate of ending up with low quality results or with low efficiencies. To obtain this, it may often be a good solution that search structures characterized by large and highly overlapped subspaces are substituted by leaner, less overlapped, recursively recombined search ones, or are modified into them. In that framework, the indicators CCF, FSW, NSM, URR and IRUIF may result useful as aids for devising efficient search structures.

The most advantageous situation for search quality was found when CCF, FSW and NSM are high. The most advantageous situation for search efficiency was instead found when net search sizes are low and both FSWs and NSM are high. Which implies that a good way to seek for search efficiency may be that of looking for comparatively high CCFs, FSWs and NSMs (maybe compounded into IRUIF) in combination with small net search sizes. Small search sizes will probably imply that a search structure is not highly integrated. But the highest integration at equal search size should be chosen, for efficiency.

In Figure 12, Search Qualities, Search Sizes and IRUIFs are shown for all the tested recursive searches. It can there be noted on a wider sample that search quality does not only depend on search size, but vary with the kind of recursive recombination.

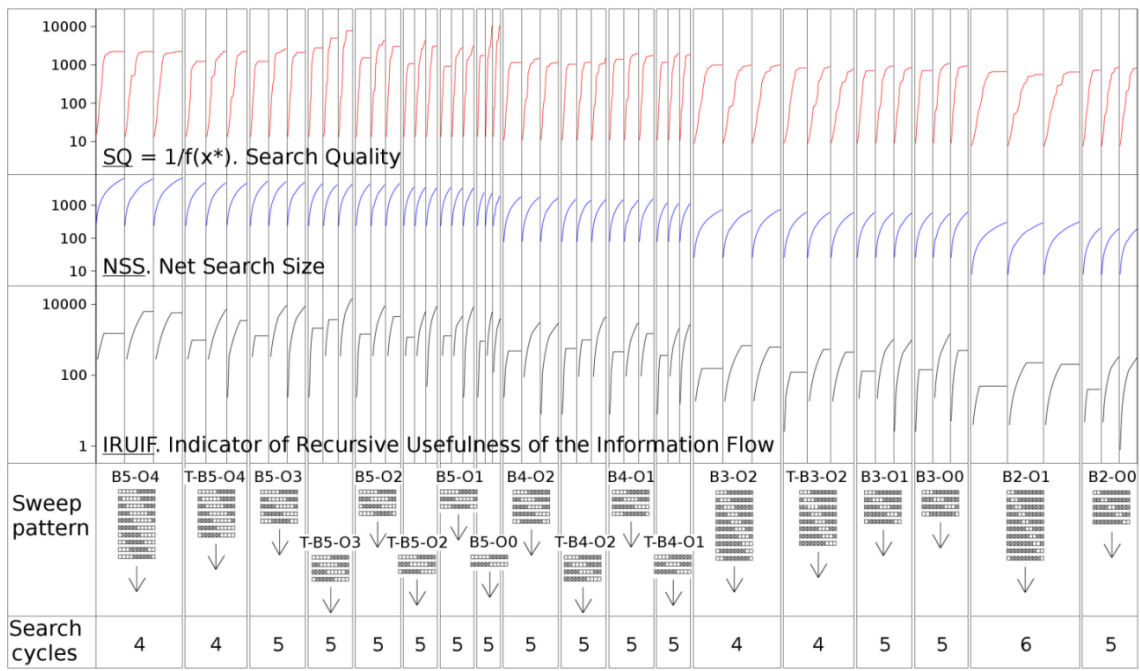

Fig. 11. Full results of the searches shown in Figures 4 and 5 regarding the recursive multi-cycle searches. The minimization problem was transformed in a maximization one. The scales are logarithmic. For each search structure, three recombination types were tested: absent (left - "n"); mild (centre - "A"); stronger (right - "B").

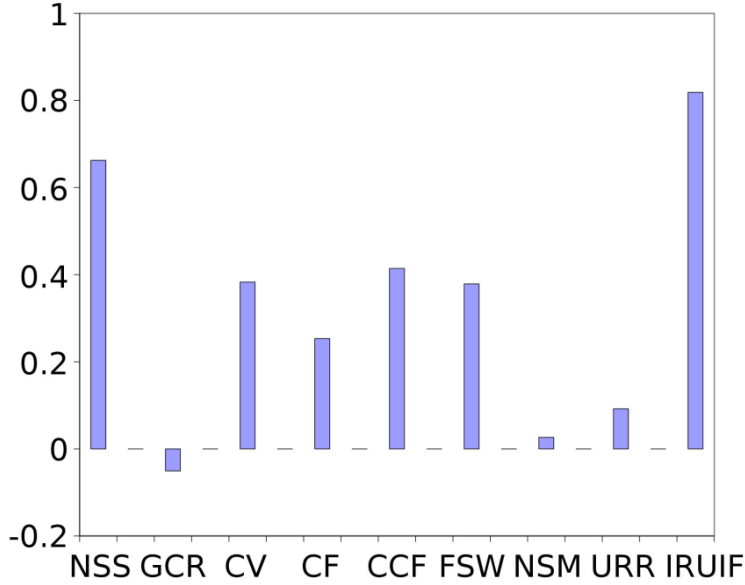

Fig. 12. Correlation between some of the functions regarding multi-cycle, recursive searches and Search Quality ($1/f(x^*)$) relative to all the trials. Legend: NSS: Net Search Size; GCR: Global Commonality Ratio; CV: Commonality Volume; CF: Commonality Flow; CCF: Cumulative Commonality Flow; FSW: Freshness of the Search Wake; NSM: Novelty of the Search Move; URR: Usefulness of the Recursive Recombination, IRUIF: Indicator of Recursive Usefulness of the Information Flow.

**Table 1 - Summary of functions and attributed terminology**
The following indicators are referred to algorithm (2). The iterations taken into account are k = 1, …, K.

*Search Quality* at iteration k
$$SQ^k = \frac{1}{f(x^*)},$$
where *f(x*)* is the optimum objective function value at iteration k. This is in the case that the minimization problem is turned into a maximization one. When the problem remains a minimization one, Search Quality is:
$$SQ^k = f(x^*).$$

*Local Overlap* between the active block at iteration k and that at (k+1)
$$\Omega^{(k,k+1)} = A^{(k)} \cap A^{(k+1)}.$$

*Local Overlap Size* between the active blocks at iteration k and at (k+1)
$$LOS^{(k,k+1)} = \prod_{j=1}^{\Omega^{(k,k+1)}} K$$
Local overlap size is the size of the search happening in the overlaps between an active block at iteration k and the active block at the next iteration, (k + 1).

*Total Overlap Size* at iteration k
$$TOS^K = \prod_{k=1}^{K} LOS^{(k-1,k)}.$$

Active Block Size at iteration k
$$ABS^{(k)} = \prod_{j=1}^{A^k} K.$$

Net Search Size between blocks k and (k+z)
$$NSS^{(k,k+z)} = \prod_{n=k}^{(k+z)} ABS^{(n)}.$$

*Gross Search Size* at iteration k
$$GSS^K = \prod_{k=1}^{K} ABS^k .$$

*Net Search Size* at iteration k
$$NSS^k = GSS^k - TOS^k .$$

*Search Efficiency* at iteration k
$$SE^k = \frac{SQ^k}{NSS^k} ,$$
when the minimization problem is turned into a maximization one.

*Global Commonality Ratio* at iteration k
$$GCR^k = \frac{TOS^k}{NSS^k} .$$

*Commonality Volume* at iteration k
$$CV^k = NSS^k \, GCR^k .$$

*Local Commonality Ratio* between iterations k and (k+1)
$$LCR^{(k,k+1)} = \frac{SS^{(k,k+1)}}{LOS^{(k,k+1)}} .$$

*Commonality Flow* at iteration k
$$CF^k = \sqrt[3]{\prod_{n=1}^{k} CF^k \, LCR^k \, ABS^k} .$$

*Cumulative Commonality Flow* at iteration k
$$CCF^k = \sum_{j=1}^{K} CF^k .$$

*Sum of the Ages in the Search Wake* at iteration k
$$SASW^k = \sum_{j=1}^{m} a_{x_j} ,$$
where $a_{x_j}$ is the age of parameter $x_j$ of $x_1, \ldots, x_m$, in terms of number of subspace iterations. At iteration k, if $x_j$ has been updated the last time at iteration t, z = k – t.

*Average Age of the Search Wake* at iteration k
$$AASW^k = \frac{SASW_k}{m} ,$$
where *m* is the number of parameters.

*Freshness of the Search Wake* at iteration k
$$FSW^k = \frac{1}{AASW^k} .$$

*Overlap Between Search Moves*: between search move from iterations k-2 to k-1 and search move from iterations k-to k. Each of the two search states (before, i.e. k-2, and after, i.e. k-1) is constituted by a list in numerical order.
$$OBSM^{(k-1,k)} = M_{(A^{(k-2,k-1)})} \cap M_{(A^{(k-1,k)})} ,$$
where $M_{(A^{(k-2,k-1)})}$ is the search move between iterations k-2 and k-1, and $M_{(A^{(k-1,k)})}$ is the search move between iterations k-1 an k.

*Elements of Novelty* in the passage from search move between iteration k-2 and k-1 to search move between iteration k-1 and k, measured in pairs constituted by elements in the same positions in the two subspaces - represented as ordered (in this specific case, numerically ascending) lists - taken from the leftmost to the rightmost.
$$EN^{(k-1,k)} = \frac{A^{(k-2)} + A^{(k-1)} + A^{(k)}}{3} - OBSM^{(k-1,k)} ,$$
where $A^{(k-2)}$ is the active block at iteration k-2, $A^{(k-1)}$ is the active block at iteration k-1, and $A^{(k)}$ is the active block block at iteration k.

*Most Novel Search Move* at iteration k
$$MNSM^k = min\ (EN^{(2,k)}, EN^{(3,k)}, \ldots, EN^{(k-1,k)})\ .$$

*Novelty of the Search Move* at iteration k
$$NSM^k = \frac{MNSM^k}{m}\ .$$

*Usefulness of Recursive Recombination* at iteration k
$$URR^k = FSW^k\ NSM^k\ .$$

Alternative definition of *Usefulness of Recursive Recombination* at iteration k
$$URR^k = \sqrt{(K)^{FSW^k}(K)^{NSM^k}}\ .$$

This version of the indicator, not adopted in the here presented tests, may be more appropriate for cases in which repetitions of non-new search moves show to produce improvements in convergence.

*Indicator of Recursive Usefulness of the Information Flow* at iteration k
$$IRUIF^k = CCF^k\ (URR^k)^2\ .$$

## 5. Implications of the indicators and possible applications to GAs

Both novelty of the search moves (NSM) and freshness of the search wake (FSW) are calculated on the basis of past moves. NSM makes possible to recognize the most recombination-intensive non-visited solutions in problem space; FSW informs about the rate of information renewal linked to search moves. Compounded with Cumulative Commonality Flow, they constitute a measure for evaluating the usefulness of search moves and can be used to make *designed* recombination in block-based (decomposition-based) optimization strategies better (more useful) than random one; to make, in a certain sense, artifice better than nature and guide explorations so that the strongest search moves available are always selected, by making them more integrated, fresher, and novel than random. A recursively recombined block search in which the selection of active blocks were supported by those indicators (by choosing the active blocks from a large set of randomly generated ones) would compound the advantage of being able to perform exploration without disruption with the possibility of profiting from some memory of the past. The here presented experiments suggest that this strategy would produce better performances; likely, by improving their ability to deal with discontinuities and non-linearities.

The limitation of this strategy is that its advantage would likely only be present when exploring new portions of problem space. Which always happens with block searches

(due to their inability to backtrack); but may not happen with EAs, due to their ability to "population-backtrack". The advantage that those indicators may bring to EAs may therefore specifically regard the exploration phases in which novel combinations are tried. (It is worth noting that the adoption of different algorithmic strategies for different computation phases is an already practiced solution in GAs - De Jong, 1992.)

**6. Conclusive remarks**

Model mutation and decomposition-based search may contribute to the settling of some form of structured and automated of "dwelling on problems" in building and urban design. Dwelling on problems implies iteration of analyses and recombination of perspectives. The two form a powerful combination to pursue quality efficiently, but require time. Automated processes condense the expenditure of time. The availability of some form of automated "dwelling on problems" may inform designers of the consequences of their options even before they begin to think. On the theoretical side, this would imply a leap of abstraction into a meta-approach in which an essential part of creative design is constituted by thinking to thinking to it.

In the presented decomposition-based computational experiments, recombination resulted advantageous for the pursuit of both search quality and efficiency, because it produced commonality in a dynamic manner, which eases convergence. Heuristic indicators were defined to support the design of recursive block searches, linking search structure integration, amount of refreshment of the search wake, amount of recombination between search moves, and search size, showing a good correlation with search quality and efficiency. Those indicators showed that the conditions for a recursive search to yield high quality results are that the search structure is well integrated, the search wake is fresh and the search moves are novel; and that the conditions for a recursive search to be efficient are that the net search size is small, the search wake is fresh and the search moves are novel.

Results also suggest that the search performances of "well-designed" search structures are less dependent than other ones on subspace decomposition; specifically, from criteria of assignment of parameters to subspaces. Which is likely to be useful, because it can make them more robust, by making their design less a delicate operation. And which may have useful implications on problem decomposition as well, meaning that, all other things being equal, recursive recombination allows smaller overlaps between subspaces, and therefore smaller search sizes, to obtain a same amount of search efficiency and quality.

Because the proposed indicators measure recombination and search path exploitation, they may be used in decomposition-based optimization searches to select the most advantageous search moves, in the sense of the most integrated, freshest, and most novel. Deriving from this possibility is the practical strategy of combining design problems in reasonably small, overlapping information units, which may be devised by (1) mapping the problems the way they are represented through abstraction, and/or following reasons deriving from (2) multidisciplinary optimization (MDO), (3) decomposition by resolution, and (4) decomposition by system; then visiting them though recursively recombined subspace searches, to seek for efficiency. This strategy may allow to reduce exponential search growth when dealing with large and/or complex search structures, possibly composed by hierarchically organized problems of problems.

Such a strategy may also be useful to support a bottom-up approach to design (Graham, 1993), because it makes possible that active subspaces (i.e. subproblems: model features or performance aspects) in the first phases of design begin to be explored one after the other even before the following ones have been planned; which compounds the possibility of an organic growth of a project and of the inquiries supporting it. That approach may also make possible for designers to reach for search efficiency without abiding their procedural habits.

It is likely that in recursive block searches the adoption of (even population-based) global search strategies used at local level different from linear search and local search strategies linking global subspace searches different from coordinate search could produce incremental efficiency improvements. It would next be worth verifying if - and in what measure - the defined indicators are suitable for making operations like mutation and recombination more advantageous than random ones even for evolutionary algorithms.

It may be finally worth noting that the conclusions reached in the present research may not be considered unrelated to human-decision-driven design. Because recursive decomposition-based methods are not population-based, they are likely to be more suitable to map (represent) human design processes than genetic algorithms (which, on the other hand may be more suitable to model cultural facts). This is the approach taken by rational analysis in cognitive modelling (Anderson, 1991). Seen from that perspective, the study of indicators like the ones here defined may even appear of some use for cultivating an intuition for the probability of innovation potential linked to design choices. In that sense, the obtained results may be interpreted as clues showing the cause of the evidence that bad design is likely to take place when information is poorly exchanged and recombined among agents (cores of reasoning about decisional matters, actors in decisional processes, or decisional entities involving actors) whilst good design is most likely to happen through a frantically recombined information exchange process (in which actors share their results before they are fully exploited and keep poking into each other's business), which maximizes the potential of the intellectual resources in place. And they may as well be interpreted as clues about why, when a certain design problem can be fit into the mind of a single person or small group, it is not rare that the single person or the small group outperforms large groups or committees, in spite of the smaller sheer computing power and search volume. Recursive recombination allows a design search process to escape the fate written in its size: in the former case, information can be more thoroughly recombined and integrated.

## 7. Conventions and symbols

Vectors are always column vectors, and their elements are denoted by superscripts. Elements of a set or a sequence are denoted by subscripts.

| | |
|---|---|
| $a \in A$ | a is an element of A. |
| $A \subset B$ | A is a subset of B. |
| $A \cap B$ | intersection of the sets A and B. |
| $A \setminus B$ | Set A minus set B. |
| $x^*$ | optimum objective function vector. |

## 8. Acknowledgements

I would like to thank Prof. Gianni Scudo, Politecnico di Milano, for accepting the use of an intermediate versions of my software tools in the context a research of which he was coordinator during 2011 and 2012: "Systemic integration of renewable energy sources technologies in the built environment", funded by the Italian Ministry of University and Research, which has preceded the development of the present research. I would also like to thank my colleague Dr.Eng. Elisabetta Rotta for reviewing the formulas in Table 1. Finally, very important, thanks to some blind reviewer around the world.

## 9. Appendix 1
**Further information on the adopted modelling and exploration criteria**

Since the model mutations were made at constant volume, a decrease of building width or depth corresponded to an increase of building height (Fig. 1). To obtain warping at constant height, side walls were rotated in plan with mirror symmetry, each with respect to its centre point. A negative warping degree in Tables 2 and 3 in Supplementary Materials indicates a "convexity" in the direction perpendicular to the back wall in plan; and negative one, the opposite. Orientation is expressed as rotation from south, negative when counterclockwise.

Three possible envelope construction solutions (the same for walls, floor and ceiling) were taken into account, having the same thermal transmittance and solar absorptance (0.1, due to a white render layer on both faces). The three types were derived by varying the order of the massive layers, to obtain different effective capacities. The layers of the "heavy" envelope were, from inside out, 12 cm brick, 12 cm brick, 10 cm rock wool. Those of the "medium" one were 12 cm brick, 10 cm rock wool, 12 cm brick. Those of the "light" one were 10 cm rock wool, 12 cm brick, 12 cm brick. Further information is made available on the web (Brunetti, 2014).

In winter, the building was heated with an ideal system in which thermal transmission was 80% convective and 20% radiant, activated over 18 °C. In summer, it run free-floating. A mass-flow network was present, integrating windows which opened above 26 °C. Solar shading devices for glazed openings were taken into account through blind/shutter controls also activated over 26 °C. The rules that were used to update pressure coefficients on walls were conceived more to obtain that the movements of the model and the changes in building's height clearly affected wind exposure of windows than to set up realistic simulations conditions. The conditions were that if a facade obstruction was higher than the building, and if the distance between the upper edge of the obstruction and the upper edges of the building was less than two times the height of the obstacle, the facade was considered partially exposed to wind; and if that distance was less than 1, it was considered obstructed. If those conditions were not met, the facade was considered fully exposed to wind.

Figure 10 in Supplementary Materials shows the distribution of performances relative to each parameter. The combination of parameters that obtained the best and worst performances are shown in Tables 2 and 3 in Supplementary Materials. The positions x and y on site could likely be causes of entrapment into local optima, due to obstructions. Reductions in building width and depth and the "light", low-mass construction solution produced low performances.

| Building width, variation, m | Building depth, variation, m | Warping angle, variation, ° | Front window, % of wall | Side windows, % of wall | Back window, % of wall | South-north position relative to centre, variation, m | East-west position relative to centre, variation, m | Rotation relative to south, variation, ° | Constr. solution | February heating load, kWh | July average max resultant temp., °C | Objective function |
|---|---|---|---|---|---|---|---|---|---|---|---|---|
| 0 | 0 | +36 | 25 | 5 | 5 | -9 | +6 | 0 | heavy | 380.2 | 30.34 | 0 |
| 0 | 0 | -36 | 25 | 5 | 5 | -9 | +6 | 0 | heavy | 331.8 | 31.15 | 0.1020 |
| 0 | 0 | +36 | 25 | 5 | 5 | 0 | +6 | 0 | heavy | 393.2 | 30.25 | 0.1025 |
| 0 | 0 | +36 | 25 | 5 | 5 | -9 | +6 | -45 | heavy | 393.9 | 30.3 | 0.1639 |
| 0 | 0 | -36 | 25 | 5 | 10 | -9 | +6 | 0 | heavy | 338 | 31.13 | 0.1742 |
| 0 | 0 | +36 | 15 | 5 | 5 | -9 | +6 | 0 | heavy | 409.9 | 30.08 | 0.1789 |
| 0 | 0 | -36 | 25 | 5 | 5 | -9 | +6 | 0 | heavy | 356.8 | 30.86 | 0.1801 |
| 0 | 0 | 0 | 25 | 5 | 5 | -9 | +6 | 0 | medium | 358.7 | 30.84 | 0.1881 |
| 0 | 0 | +36 | 5 | 5 | 5 | -9 | +6 | 0 | heavy | 440.1 | 29.65 | 0.1921 |
| 0 | 0 | -36 | 15 | 5 | 5 | -9 | +6 | 0 | heavy | 377.5 | 30.57 | 0.1939 |

Table 2. Top 10 results of the linear exploration regarding the 10 considered design parameters. The objective function was here to be minimized.

| Building width, variation, m | Building depth, variation, m | Warping angle, variation, ° | Front window, % of wall | Side windows, % of wall | Back window, % of wall | South-north position relative to centre, variation, m | East-west position relative to centre, variation, m | Rotation relative to south, variation, ° | Constr. solution | February heating load, kWh | July average max resultant temp., °C | Objective function |
|---|---|---|---|---|---|---|---|---|---|---|---|---|
| -2 | -2 | +36 | 25 | 15 | 15 | 0 | +6 | +45 | light | 670.3 | 35.23 | 0.9317 |
| -2 | -2 | 0 | 25 | 15 | 15 | +9 | +6 | +45 | light | 719.8 | 34.54 | 0.9353 |
| -2 | -2 | -36 | 25 | 15 | 15 | +9 | 0 | +45 | light | 729.7 | 34.4 | 0.9359 |
| -2 | -2 | -36 | 25 | 15 | 15 | +9 | -6 | 0 | light | 748.9 | 34.14 | 0.9381 |
| -2 | -2 | -36 | 25 | 15 | 15 | +9 | +6 | +45 | light | 708.1 | 34.74 | 0.9382 |
| -2 | -2 | -36 | 25 | 15 | 10 | +9 | -6 | +45 | light | 741.9 | 34.25 | 0.9388 |
| -2 | -2 | -36 | 25 | 15 | 15 | 0 | -6 | +45 | light | 722 | 34.63 | 0.9478 |
| -2 | -2 | 0 | 25 | 15 | 15 | +9 | -6 | +45 | light | 764.4 | 34.1 | 0.9572 |
| -2 | -2 | +36 | 25 | 15 | 15 | +9 | -6 | +45 | light | 747.8 | 34.76 | 0.9996 |
| -2 | -2 | -36 | 25 | 15 | 15 | +9 | -6 | +45 | light | 764.4 | 34.52 | 1 |

Table 3. Bottom 10 results of the linear exploration regarding the 10 considered design parameters. The objective function was here to be minimized.